\documentclass[letterpaper, 10 pt, conference]{ieeeconf}  
\usepackage[dvipsnames]{xcolor}

\IEEEoverridecommandlockouts                              

\usepackage{amsmath}
\usepackage{amssymb}
\usepackage{amsfonts}

\usepackage{graphicx}

\title{\LARGE \bf
PH-KAN: Port-Hamiltonian Kolmogorov-Arnold Network$^{*}$
}

\author{Achraf El Messaoudi$^{1}$, Karim Cherifi$^{1}$, Yann Le Gorrec$^{1}$, Yongxin Wu$^{1}$
\thanks{$^{1}$ Université Marie et Louis Pasteur, SUPMICROTECH, CNRS, institut FEMTO-ST, F-25000 Besançon, France}
\thanks{$^{*}$ This work has been supported by the
EIPHI Graduate School (contract: ANR-17-EURE-0002) and Bourgogne-Franche-Comté Region.}
}

\begin{document}

\maketitle
\thispagestyle{empty}
\pagestyle{empty}

\begin{abstract}
Data-driven machine learning approaches have become increasingly attractive for nonlinear system identification, but standard models often fail to preserve the underlying physical structure and remain difficult to interpret, especially when no analytical model is available. In this context, port-Hamiltonian (pH) models provide a natural physics-informed representation. However, when these models are parameterized with standard multilayer perceptrons (MLPs), the learned constitutive components often remain poorly interpretable. In this paper, we propose a structure-preserving identification framework for nonlinear port-Hamiltonian systems based on Kolmogorov-Arnold Networks (KANs). The proposed PH-KAN model parameterizes the interconnection matrix, dissipation matrix, Hamiltonian, and input mapping using dedicated KAN blocks, while enforcing the port-Hamiltonian constraints by construction. This yields constitutive representations in which the nonlinear functions defining the identified pH components can be explicitly inspected, leading to a more interpretable model than with standard MLP-based parameterizations.
\end{abstract}

\section{INTRODUCTION}

Nonlinear system identification is a fundamental problem in control and engineering, especially in situations where no analytical model is known a priori and the dynamical model must therefore be inferred directly from available data. In this context, data-driven machine learning methods have become increasingly attractive due to their ability to approximate complex nonlinear dynamics with high accuracy \cite{pillonetto2025}. However, standard black-box models often fail to preserve the underlying physical structure and remain difficult to interpret, which limits their use in applications where physical consistency and model transparency are important.

For many physical systems, an energy-based representation provides a natural modeling framework. In recent years, structure-preserving neural approaches such as Hamiltonian and Lagrangian neural networks have highlighted the benefits of embedding physical priors into learning-based dynamical models \cite{greydanus2019,cranmer2020}. However, these approaches mainly target conservative systems and do not explicitly account for dissipation and external inputs, which are essential in many engineering applications.

In this context, the port--Hamiltonian (pH) framework provides a particularly appealing representation. A pH model explicitly incorporates stored energy, dissipative effects, interconnection structure, and external inputs, making it well suited to the modeling and control of nonlinear multiphysical systems \cite{vanderschaft2014, NelsonCisneros_2025, AmalHammoud_2026}. This has motivated the development of port-Hamiltonian neural networks, which aim to identify pH-consistent models from data while preserving the structural constraints of the formalism \cite{desai2021,moradi2026,cherifi2025}. In such approaches, the goal is not only to fit an observed vector field, but also to recover constitutive objects such as the Hamiltonian, the dissipation term, and the interconnection structure in a physically meaningful way.

Despite these advantages, most existing neural approaches to pH system identification rely on multilayer perceptrons (MLPs) to model the nonlinear pH structure. While effective, MLP-based parameterizations remain difficult to interpret, so that the learned constitutive relations often retain a black-box character even when the overall model satisfies the pH constraints. This is a significant limitation in the pH setting, where interpretability is one of the key motivations for adopting structured models in the first place.

Kolmogorov--Arnold Networks (KANs) have recently emerged as a promising alternative to MLPs \cite{liu2024}. By placing learnable univariate functions on edges rather than fixed activations on nodes, KANs offer a representation that is more intepretable. Recent work has also begun to explore KANs for nonlinear system identification \cite{cruz2025}. A closely related recent work \cite{cherifi2025} also considers KAN-based parameterizations for nonlinear port--Hamiltonian system identification and compares them with MLP-based alternatives. However, this implementation based on FastKAN \cite{li2024fastkan}, does not support the direct inspection of learned functions used here for interpretability analysis. Moreover, that work does not investigate whether the favorable scaling behavior reported in the KAN literature remains valid in a structure-preserving setting.

Motivated by these observations, we propose a structure-preserving KAN-based framework for nonlinear port--Hamiltonian system identification, aimed at obtaining explicit and interpretable representations of the learned nonlinear pH components and studying their scaling behavior under port-Hamiltonian constraints. 

The contributions of this paper are as follows:
\begin{itemize}
    \item We propose a Port--Hamiltonian Kolmogorov--Arnold Network (PH-KAN) framework for nonlinear system identification, in which the pH structure is preserved by construction while the interconnection matrix, dissipation matrix, Hamiltonian, and input mapping are parameterized by KAN blocks.
    \item We show how the resulting parameterization supports interpretability of the learned model through the inspection of edge functions and the use of sparsification and pruning mechanisms.
    \item We evaluate the proposed approach on a port--Hamiltonian benchmark and show that, in addition to its interpretability benefits, PH-KAN can outperform a structure-matched PH-MLP baseline in predictive accuracy for certain training regimes and model sizes.
\end{itemize}

\section{PROBLEM FORMULATION}

\subsection{Port-Hamiltonian Systems}

Port-Hamiltonian (pH) systems provide a physically meaningful representation for a broad class of nonlinear dynamical systems by explicitly separating power-conserving interconnections, dissipative effects, and external inputs. In this framework, the system dynamics are written as
\begin{equation}
\dot{x} = \big(J(x)-R(x)\big)\nabla H(x) + B(x)u, \quad y=B(x)^\top \nabla H(x)
\label{eq:ph_general}
\end{equation}

where $x \in \mathbb{R}^n$ denotes the state, $u \in \mathbb{R}^m$ is the input, $y \in \mathbb{R}^m$ is the power conjugated output, $H:\mathbb{R}^n \to \mathbb{R}$ is the Hamiltonian function, $J(x)\in\mathbb{R}^{n\times n}$ is the interconnection matrix, $R(x)\in\mathbb{R}^{n\times n}$ is the dissipation matrix, and $B(x)\in\mathbb{R}^{n\times m}$ is the input matrix.

The pH structure imposes the following constraints:
\begin{equation}
J(x) = -J(x)^\top, \qquad R(x)=R(x)^\top \succeq 0.
\label{eq:ph_constraints}
\end{equation}
The skew-symmetry of $J(x)$ ensures that the interconnection term redistributes energy without generating or dissipating it, while the positive semidefiniteness of $R(x)$ models irreversible energy dissipation. The Hamiltonian $H(x)$ typically represents the total stored energy of the system, and its gradient $\nabla H(x)$ plays the role of the generalized effort vector.

An important consequence of \eqref{eq:ph_general} is the energy balance relation
\begin{equation}
\dot{H}(x)
=
-\nabla H(x)^\top R(x)\nabla H(x)
+
y^\top  u \preceq y^\top  u ,
\label{eq:energy_balance}
\end{equation}
which makes explicit how the stored energy evolves due to dissipation and external power injection through the power-conjugate input and output ports. This property makes pH models particularly attractive for identification and control, since they embed physically interpretable structure directly into the dynamics.

In this work, we focus on nonlinear systems that admit a pH representation of the form \eqref{eq:ph_general} and seek to identify the constitutive objects $J(x)$, $R(x)$, $H(x)$, and $B(x)$ from data while preserving the pH structural properties by construction.

\subsection{Identification Problem Setting}

Consider a nonlinear control system for which input-state trajectory data are available. Given $N$ samples
\begin{equation}
\mathcal{D}=\{(x_k,u_k,\dot{x}_k)\}_{k=1}^N,
\label{eq:dataset}
\end{equation}
our objective is to identify a structured model of the form
\begin{equation}
\dot{\hat{x}} = \big(\hat{J}(x)-\hat{R}(x)\big)\nabla \hat{H}(x) + \hat{B}(x)u,
\label{eq:identified_model}
\end{equation}
where $\hat{J}(x)$, $\hat{R}(x)$, $\hat{H}(x)$, and $\hat{B}(x)$ are learned from data.

The identification task is not merely to fit the observed vector field, but to recover a model that remains consistent with the pH formalism. To this end, the learned quantities must satisfy
\begin{equation}
\hat{J}(x) = -\hat{J}(x)^\top, \qquad
\hat{R}(x)=\hat{R}(x)^\top \succeq 0
\label{eq:learned_constraints}
\end{equation}
Rather than enforcing these properties only through soft penalties, our objective is to impose them directly through the model parameterization. This yields a structure-preserving identification framework in which physical consistency is embedded into the learned dynamics.

Let $\hat{f}(x,u)$ denote the identified vector field,
\begin{equation}
\hat{f}(x,u) := \big(\hat{J}(x)-\hat{R}(x)\big)\nabla \hat{H}(x) + \hat{B}(x)u.
\label{eq:identified_vector_field}
\end{equation}
The identification problem can then be stated as follows.

\medskip
\noindent\textbf{Problem:}
Given trajectory data generated by an unknown nonlinear pH system, learn structured approximations $\hat{J}(x)$, $\hat{R}(x)$, $\hat{H}(x)$, and $\hat{B}(x)$ such that:
\begin{enumerate}
    \item the predicted state derivative $\hat{f}(x,u)$ matches the observed dynamics accurately.
    \item the learned model preserves the key pH structural constraints.
    \item the resulting representation is adequate for interpretation.
\end{enumerate}

To instantiate this structured approximation problem, we consider Kolmogorov-Arnold Networks (KANs) as nonlinear approximators for the unknown components of the model. KANs are introduced here as an alternative to multilayer perceptrons (MLPs), which are widely used in the literature for neural dynamical system identification and structured neural modeling. Beyond their approximation capabilities, KANs are of particular interest because they offer the potential for improved interpretability compared with standard MLPs, which typically behave as more opaque black-box models. In addition, recent KAN literature reports favorable scaling behavior relative to MLPs, which motivates investigating their use in the present structure-preserving identification setting. 

\subsection{PMSM Benchmark}

To assess the proposed identification framework, we consider a Permanent Magnet Synchronous Motor (PMSM) described in port-Hamiltonian form \cite{L2gain_book}.
The state is defined as
\begin{equation}
x = \begin{bmatrix}\phi_d & \phi_q & p\end{bmatrix}^\top,
\end{equation}
where $\phi_d$ and $\phi_q$ are the $d$- and $q$-axis flux linkages, and $p$ is the angular momentum. The control input is
\begin{equation}
u = \begin{bmatrix}u_1 & u_2\end{bmatrix}^\top,
\end{equation}
corresponding to the electrical inputs applied along the two axes.

The PMSM dynamics are given by
\begin{equation}
\dot{x} = \big(J(x)-R\big)\nabla H(x) + Bu,
\label{eq:pmsm_ph}
\end{equation}
with Hamiltonian
\begin{equation}
H(x)=\frac{1}{2L}\left(\phi_d^2+\phi_q^2\right)+\frac{1}{2J_m}p^2,
\label{eq:pmsm_H}
\end{equation}
and
\begin{equation}
\begin{aligned}
J(x)&=
\begin{bmatrix}
0 & 0 & \phi_q\\
0 & 0 & -(\phi_d+\Phi)\\
-\phi_q & \phi_d+\Phi & 0
\end{bmatrix},
\;\;
R=
\begin{bmatrix}
r & 0 & 0\\
0 & r & 0\\
0 & 0 & \beta
\end{bmatrix}\\
B&=
\begin{bmatrix}
1 & 0\\
0 & 1\\
0 & 0
\end{bmatrix},
\qquad
\nabla H(x)=
\begin{bmatrix}
\phi_d/L\\
\phi_q/L\\
p/J_m
\end{bmatrix}.
\end{aligned}
\label{eq:pmsm_matrices}
\end{equation}
The physical parameters used in this work are
\begin{equation}
\begin{aligned}
J_m &= 0.012,&
L &= 3.8\times10^{-3},&
\beta &= 0.0026,\\
r &= 0.225,&
\Phi &= 0.17. & &
\end{aligned}
\label{eq:pmsm_params}
\end{equation}
Here, $J_m$ denotes the rotor inertia, $L$ the phase inductance, $\beta$ the viscous friction coefficient, $r$ the phase resistance, and $\Phi$ the constant rotor magnetic flux.

The PMSM constitutes a representative benchmark for evaluating the proposed identification framework as it is widely used in high-performance motion control applications, including robotics, electric drives, and electric vehicles, and their dynamics are well understood in the control literature. In addition, an interpretable representation is highly desirable for PMSM control, as many advanced control strategies, including passivity-based control and energy-shaping approaches, explicitly rely on the underlying physical structure of the model. Learning an interpretable port-Hamiltonian form therefore enables the integration of data-driven modeling with structure-preserving control design while maintaining physical consistency and facilitating analysis of the learned dynamics.

In the following section, we introduce the proposed PH-KAN framework for directly identifying the components of this model from trajectory data.

\section{Proposed PH-KAN Identification Framework}

\subsection{Kolmogorov-Arnold Networks}
A Kolmogorov--Arnold Network (KAN) is a nonlinear map obtained by composing successive layers. Let
$
x \in \mathbb{R}^{n_0}
$
denote the input vector. A KAN with $L$ layers defines a function
\begin{equation}
\mathrm{KAN}(x)
=
\left(\Phi_{L-1}\circ \Phi_{L-2}\circ \cdots \circ \Phi_0\right)(x),
\label{eq:kan_composition}
\end{equation}
where each layer
\[
\Phi_\ell : \mathbb{R}^{n_\ell} \to \mathbb{R}^{n_{\ell+1}}
\]
maps an input vector $z^{(\ell)} \in \mathbb{R}^{n_\ell}$ to an output vector $z^{(\ell+1)} \in \mathbb{R}^{n_{\ell+1}}$, with $z^{(0)}=x$. The symbol $\circ$ denotes function composition.

In contrast to a MLP, where scalar nonlinear activation functions are attached to the nodes, a KAN places learnable scalar nonlinearities on the edges between nodes. More precisely, for a given layer $\Phi_\ell$ and input
\[
z^{(\ell)}=\big(z^{(\ell)}_1,\dots,z^{(\ell)}_{n_\ell}\big)^\top,
\]
the $j$-th output component is defined as
\begin{equation}
z^{(\ell+1)}_j
=
\sum_{i=1}^{n_\ell}
\phi_{\ell,j,i}\!\left(z^{(\ell)}_i\right),
\qquad j=1,\dots,n_{\ell+1},
\label{eq:kan_layer}
\end{equation}
where each
\[
\phi_{\ell,j,i}:\mathbb{R}\to\mathbb{R}
\]
is a learnable univariate function associated with the edge connecting input node $i$ of layer $\ell$ to output node $j$ of layer $\ell+1$. Figure~1 illustrates this construction on a simple example: the boxes placed on the edges represent the scalar functions $\phi_{\ell,j,i}$, while the labels $\Phi_1$ and $\Phi_2$ denote the corresponding layer maps obtained by summing all edge contributions according to \eqref{eq:kan_layer}.

In the \texttt{pykan} implementation \cite{liu2024} used in this work, each edge function is written as the sum of a residual base function and a spline function,
\begin{equation}
\phi_{\ell,j,i}(\xi)
=
w_{\ell,j,i}
\left(
b(\xi)
+
\sum_{r} c_{\ell,j,i,r} B_r(\xi)
\right),
\label{eq:kan_edge_function}
\end{equation}
where $\xi\in\mathbb{R}$ is the scalar input of the edge function, $b(\cdot)$ is a fixed base function, and $B_r(\cdot)$ are one-dimensional B-spline basis functions with trainable coefficients $c_{\ell,j,i,r}$. The base function is chosen as the SiLU nonlinearity. The spline basis is defined on a one-dimensional grid, whose resolution is controlled by the grid size.

This parameterization makes each edge function flexible while preserving the layered KAN structure of \eqref{eq:kan_layer}. In practice, the main architectural hyperparameters are the layer widths $[n_0,n_1,\dots,n_L]$, which determine the network depth and width, and the grid size used in the spline expansion. In the present work, these correspond to the \texttt{width} and \texttt{grid} parameters of \texttt{pykan}. Dedicated KAN blocks are then used to parameterize the constitutive components of the identified port-Hamiltonian model.

\begin{figure*}[t]
    \centering
    \includegraphics[width=0.8\textwidth]{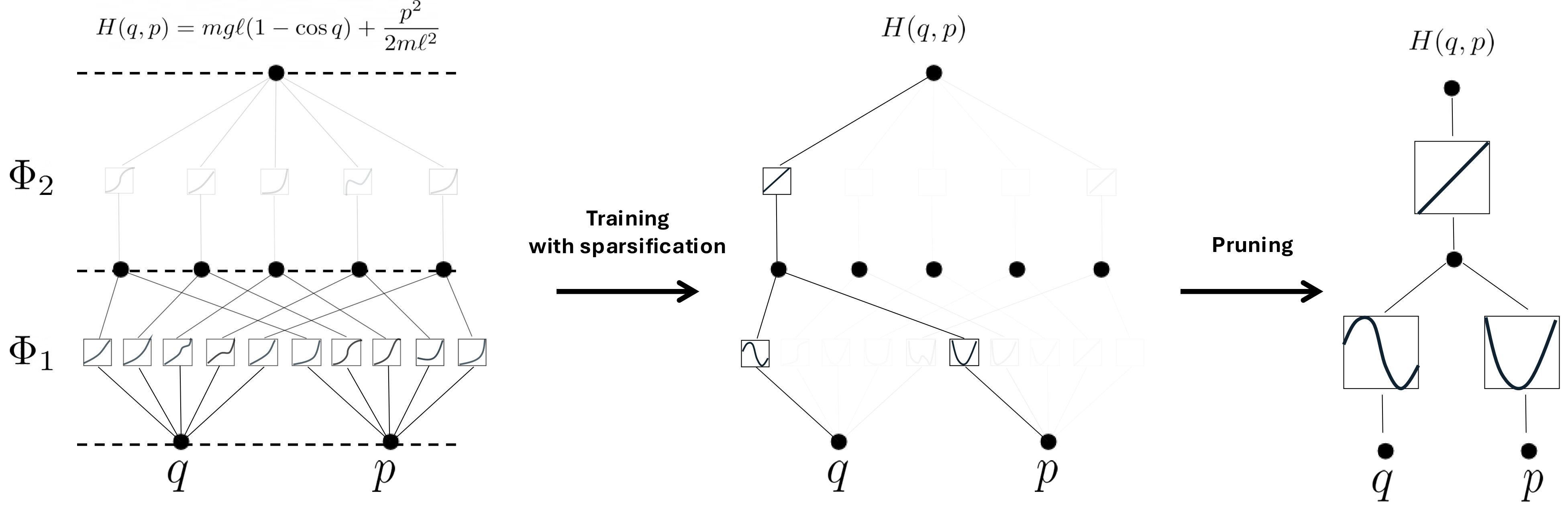}
    \caption{Illustrative example of sparsification and pruning on a KAN approximating the Hamiltonian of a simple pendulum. Weakly contributing connections are progressively eliminated, leading to a compact and more interpretable subnetwork.}
    \label{fig:kan_pruning}
\end{figure*}

\subsection{Structure-Preserving Parameterization of $J$, $R$, $H$, and $B$}

In order to enforce the pH structure, we parametrize $J(x)$, $R(x)$, $H(x)$, and $B(x)$ in a way that naturally enforces the structure.

Let
\[
M_J(x) := \operatorname{mat}_{n,n}(a_J(x)), \qquad
M_R(x) := \operatorname{mat}_{n,n}(a_R(x)),
\]
\[
h_{\mathrm{raw}}(x) := a_H(x), \qquad
b_B(x) := a_B(x),
\]
where $\operatorname{mat}_{p,q}(\cdot)$ reshapes a vector into a $p\times q$ matrix.

The interconnection matrix is parameterized as
\begin{equation}
\hat{J}(x) = M_J(x) - M_J(x)^\top,
\label{eq:J_param_compact}
\end{equation}
which guarantees
\[
\hat{J}(x) = -\hat{J}(x)^\top.
\]

The dissipation matrix is parameterized as
\begin{equation}
\hat{R}(x) = M_R(x)M_R(x)^\top,
\label{eq:R_param_compact}
\end{equation}
which guarantees
\[
\hat{R}(x)=\hat{R}(x)^\top \succeq 0.
\]

To ensure cyclo-passivity, the learned Hamiltonian is obtained from the raw scalar output $h_{\mathrm{raw}}(x)$ of the KAN block through a shifted ELU transformation:
\begin{equation}
\hat{H}(x)=\operatorname{ELU}\!\left(h_{\mathrm{raw}}(x)-b\right)+b,
\label{eq:H_elu_compact}
\end{equation}
where $h_{\mathrm{raw}}(x)$ denotes the unconstrained scalar output of the KAN block associated with $H$, and
\begin{equation}
b = H_{\min}+1.
\label{eq:b_shift_choice}
\end{equation}
Here, $\operatorname{ELU}(\cdot)$ denotes the exponential linear unit with parameter $\alpha=1$, i.e.,
\begin{equation}
\operatorname{ELU}(z)=
\begin{cases}
z, & z>0,\\
e^z-1, & z\leq 0.
\end{cases}
\label{eq:elu_def}
\end{equation}
Since $\operatorname{ELU}(z)\geq -1$ for all $z\in\mathbb{R}$, it follows directly from \eqref{eq:H_elu_compact} that
\[
\hat{H}(x)\geq -1+b = H_{\min},
\]
so that the learned Hamiltonian is bounded from below by construction.

The Hamiltonian gradient is then obtained by automatic differentiation:
\begin{equation}
\nabla \hat{H}(x)=\frac{\partial \hat{H}(x)}{\partial x}.
\label{eq:gradH_compact}
\end{equation}

Finally, the input matrix is defined as
\begin{equation}
\hat{B}(x)=\operatorname{mat}_{n,m}(b_B(x)),
\label{eq:B_param_compact}
\end{equation}
so that its dimension is directly determined by the number of states $n$ and inputs $m$.

The resulting identified model is
\begin{equation}
\dot{\hat{x}}=
\big(\hat{J}(x)-\hat{R}(x)\big)\nabla \hat{H}(x)+\hat{B}(x)u.
\label{eq:phkan_final_model}
\end{equation}

\subsection{Training Objective and Regularization}

Given a training dataset
\[
\mathcal{D}=\{(x_k,u_k,\dot{x}_k)\}_{k=1}^N,
\]
the PH-KAN model is trained by minimizing a derivative fitting loss of the form
\begin{equation}
\mathcal{L}_{\mathrm{data}}
=
\frac{1}{N}\sum_{k=1}^N
\left\|
\hat{f}(x_k,u_k)-\dot{x}_k
\right\|_2^2,
\label{eq:data_loss}
\end{equation}
where
\begin{equation}
\hat{f}(x,u)
=
\big(\hat{J}(x)-\hat{R}(x)\big)\nabla \hat{H}(x)
+
\hat{B}(x)u.
\label{eq:fhat_training}
\end{equation}

To promote simpler and more interpretable KAN representations, we add a sparsification regularization to the KAN blocks. Unlike in MLPs, where sparsity is usually imposed on linear weights, sparsity in KANs is enforced on the learnable edge functions. For a KAN edge function $\phi$, its empirical $\ell_1$ magnitude over a batch of $N_b$ samples is defined as
\begin{equation}
\|\phi\|_1
=
\frac{1}{N_b}\sum_{s=1}^{N_b}
\left|
\phi\!\left(x^{(s)}\right)
\right|.
\label{eq:phi_l1}
\end{equation}
For a KAN layer $\Phi_\ell$, the corresponding layer-wise sparsity measure is
\begin{equation}
\|\Phi_\ell\|_1
=
\sum_{i=1}^{n_\ell}\sum_{j=1}^{n_{\ell+1}}
\|\phi_{\ell,j,i}\|_1.
\label{eq:Phi_l1}
\end{equation}

As in the original KAN formulation, an entropy term is further introduced to enhance simplification. Denoting
\begin{equation}
p_{\ell,j,i}
=
\frac{\|\phi_{\ell,j,i}\|_1}{\|\Phi_\ell\|_1},
\label{eq:p_lji}
\end{equation}
the layer entropy is defined as
\begin{equation}
\mathcal{S}(\Phi_\ell)
=
-
\sum_{i=1}^{n_\ell}\sum_{j=1}^{n_{\ell+1}}
p_{\ell,j,i}\log p_{\ell,j,i}.
\label{eq:entropy_layer}
\end{equation}
In addition, we penalize the spline coefficients of the edge functions in order to favor smoother and less complex spline representations. Denoting this term by $\mathcal{C}(\Phi_\ell)$, the regularization associated with one KAN block is written as
\begin{equation}
\mathcal{R}_{\mathrm{KAN}}
=
\sum_{\ell=0}^{L-1}
\left(
\mu_1 \|\Phi_\ell\|_1
+
\mu_2 \mathcal{S}(\Phi_\ell)
+
\mu_3 \mathcal{C}(\Phi_\ell)
\right),
\label{eq:kan_regularization}
\end{equation}
where $\mu_1,\mu_2,\mu_3\geq 0$ are regularization weights.

Accordingly, the overall training loss is
\begin{equation}
\mathcal{L}
=
\mathcal{L}_{\mathrm{data}}
+
\lambda
\left(
\mathcal{R}_J
+
\mathcal{R}_R
+
\mathcal{R}_H
+
\mathcal{R}_B
\right),
\label{eq:total_loss_phkan}
\end{equation}
where $\mathcal{R}_J,\mathcal{R}_R,\mathcal{R}_H,\mathcal{R}_B$ denote the KAN regularization terms of the form $\mathcal{R}_{\mathrm{KAN}}$, associated with the parameterizations of $\hat{J}$, $\hat{R}$, $\hat{H}$, and $\hat{B}$, respectively.

After sparsification, a pruning step can be applied to extract a smaller subnetwork. The idea is to remove edge functions whose empirical contribution on the training data is negligible. For an edge function $\phi_{\ell,j,i}$, we define the empirical edge score
\begin{equation}
s_{\ell,j,i}
=
\frac{1}{N_{\mathrm{tr}}}
\sum_{s=1}^{N_{\mathrm{tr}}}
\left|
\phi_{\ell,j,i}\!\left(z^{(\ell,s)}_i\right)
\right|,
\label{eq:edge_score_pruning}
\end{equation}
where $z^{(\ell,s)}_i$ denotes the $i$-th input of layer $\ell$ for the $s$-th training sample. Given a threshold $\tau_e>0$, the edge is pruned whenever
\begin{equation}
s_{\ell,j,i}<\tau_e.
\label{eq:edge_pruning_rule}
\end{equation}

In addition, node pruning can be performed by comparing the importance of incoming and outgoing connections. A node is kept only if its incoming and outgoing scores are both above a prescribed threshold. In practice, pruning is applied after a forward pass on the training set, which provides the activation statistics required to identify weakly contributing edges and nodes. A visual illustration of this sparsification--pruning process is shown in Fig. \ref{fig:kan_pruning}.

This process yields a reduced PH-KAN model with fewer active components, improved readability, and only limited degradation in predictive accuracy. Beyond reducing the model size, pruning removes weakly contributing edges and nodes and makes the learned structure easier to inspect. In particular, it helps reveal which parts of the network remain active in the representation of the constitutive relations, thereby providing a cleaner basis for the qualitative analysis of the learned model. Optionally, the pruned model can then be retrained for a few additional optimization steps after pruning in order to recover part of the performance loss induced by the removal of weakly contributing components. 

\section{Data Generation and Experimental Setup}

\subsection{Dataset Generation}

The dataset is generated from the analytical PMSM model introduced in Section II-C. Training trajectories are obtained by numerically integrating the ground-truth dynamics
\begin{equation}
\dot{x} = \big(J(x)-R\big)\nabla H(x) + Bu
\label{eq:pmsm_data_dyn}
\end{equation}
using the RK45 integration scheme.

For each trajectory, the initial condition is sampled uniformly as
\begin{equation}
x_0 \sim \mathcal{U}\big([-0.5,0.5]\times[-0.5,0.5]\times[-1,1]\big),
\label{eq:init_cond_sampling}
\end{equation}
so as to cover a representative operating region around the origin. The two input channels are independently excited by multisine signals of the form
\begin{equation}
u_i(t)
=
a \sum_{k=1}^{N_h}
\sin\!\left(2\pi k f_0 t + \varphi_{i,k}\right),
\qquad i\in\{1,2\},
\label{eq:multisine_input}
\end{equation}
where the phases $\varphi_{i,k}$ are sampled independently and uniformly in $[0,2\pi)$.

In the experiments, we use $a=1$, $f_0=0.1$ Hz, and $N_h=40$, which yields an excitation band from $0.1$ Hz to $40$ Hz. This choice provides a persistently exciting input over a broad low-to-medium frequency range, which is well suited to the PMSM benchmark. 

We generate $200$ trajectories over a time horizon of $10$ s with sampling time $\Delta t = 0.01$ s. This results in a dataset of state samples $x_k$, input samples $u_k$, and state-derivative targets $\dot{x}_k$. All generated samples are then aggregated and randomly shuffled.

\subsection{Evaluation Protocol}

Model evaluation is carried out through trajectory rollout on an unseen test experiment. More precisely, after training, a fresh initial condition is sampled in the same range as in \eqref{eq:init_cond_sampling}, and a new pair of multisine inputs is generated using independent random phases.

The ground-truth trajectory is obtained by integrating the analytical PMSM dynamics, while the learned trajectory is obtained by integrating the identified PH-KAN model in port-Hamiltonian form:
\begin{equation}
\dot{x}
=
\big(\hat{J}(x)-\hat{R}(x)\big)\nabla \hat{H}(x)
+
\hat{B}(x)u.
\label{eq:rollout_eval}
\end{equation}

Both systems are simulated over the same horizon with the same initial condition and input signals.

The evaluation is based on a direct comparison between the ground-truth rollout and the learned rollout for the three state variables $(\phi_d,\phi_q,p)$. In addition to qualitative trajectory plots, we report the rollout mean squared error for each state:
\begin{equation}
\mathrm{MSE}_i
=
\frac{1}{T}
\sum_{k=1}^{T}
\big(\hat{x}_i(t_k)-x_i(t_k)\big)^2,
\label{eq:rollout_mse}
\end{equation}
where $x_i$ denotes the reference trajectory and $\hat{x}_i$ the corresponding KAN-based predicted trajectory of the $i$-th state component.

\section{Results}

\subsection{Rollout Performance}

\begin{figure}[t]
    \centering
    \includegraphics[width=\linewidth]{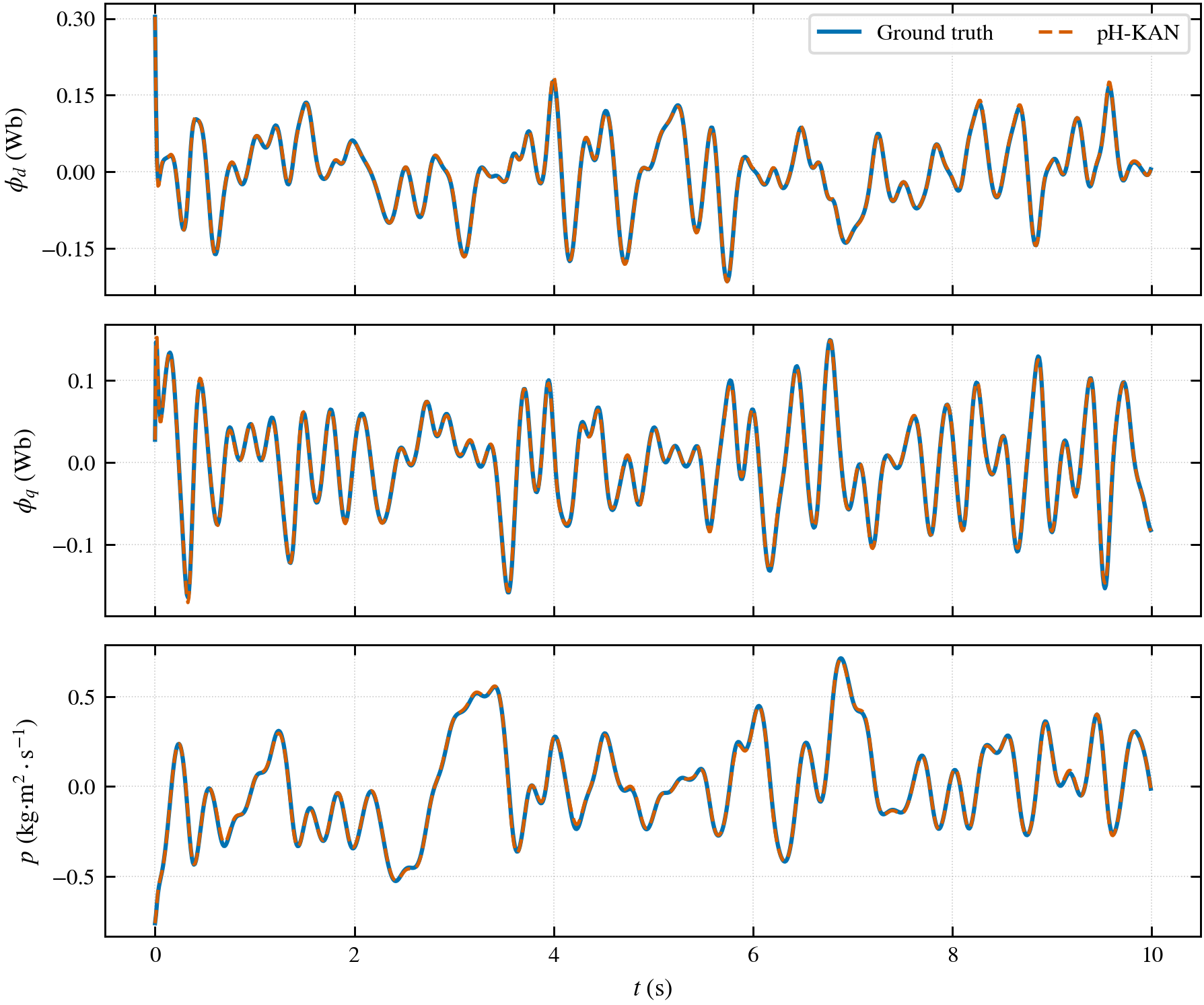}
    \caption{Comparison between the ground-truth PMSM trajectory and the PH-KAN rollout on an unseen test experiment.}
    \label{fig:rollout_comparison}
\end{figure}

Fig.~\ref{fig:rollout_comparison} shows a representative comparison between the ground-truth rollout and the PH-KAN rollout for the three state variables $(\phi_d,\phi_q,p)$. The learned model accurately reproduces the temporal evolution of all states over the full simulation horizon, indicating that it captures the coupled electromechanical behavior of the PMSM.

\begin{table}[t]
\centering
\caption{Rollout evaluation under different training configurations.}
\label{tab:rollout_configs}
\begin{tabular}{ccc}
\hline
Number of trajectories & Number of epochs & MSE \\
\hline
200 & 400 & $1.4 \times 10^{-3}$ \\
100 & 200 & $1.8 \times 10^{-2}$ \\
50  & 100 & $5.4 \times 10^{-2}$ \\
\hline
\end{tabular}
\end{table}

The different training configurations considered for rollout evaluation are summarized in Table \ref{tab:rollout_configs}. As expected, the rollout accuracy improves as the number of training samples increases, indicating that the identified PH-KAN model benefits from larger training datasets. Figure~\ref{fig:rollout_comparison} reports the rollout obtained with the configuration achieving the lowest error in Table~\ref{tab:rollout_configs}.

\subsection{Convergence Comparison with MLPs}

We now compare the parameter efficiency of the proposed PH-KAN model with that of a PH-MLP baseline. The latter uses the same port-Hamiltonian parameterization as the proposed model, i.e., the same structured decomposition into $\hat{J}(x)$, $\hat{R}(x)$, $\hat{H}(x)$, and $\hat{B}(x)$ with the same structural constraints imposed by construction. The only difference lies in the choice of nonlinear approximator: in the PH-MLP baseline, each unknown component is parameterized by a standard MLP. This provides a natural reference model, since MLP-based parameterizations are the most classical choice for nonlinear system identification and are closely related to existing neural approaches in the literature.

For comparison, the reference slopes $N^{-2}$ and $N^{-4}$ are also displayed, where $N$ denotes the number of trainable parameters. These rates are shown as illustrative neural scaling references, in line with the scaling discussion in the KAN literature \cite{liu2024,koenig2024}.

\begin{figure}[t]
    \centering
    \includegraphics[width=0.8\linewidth]{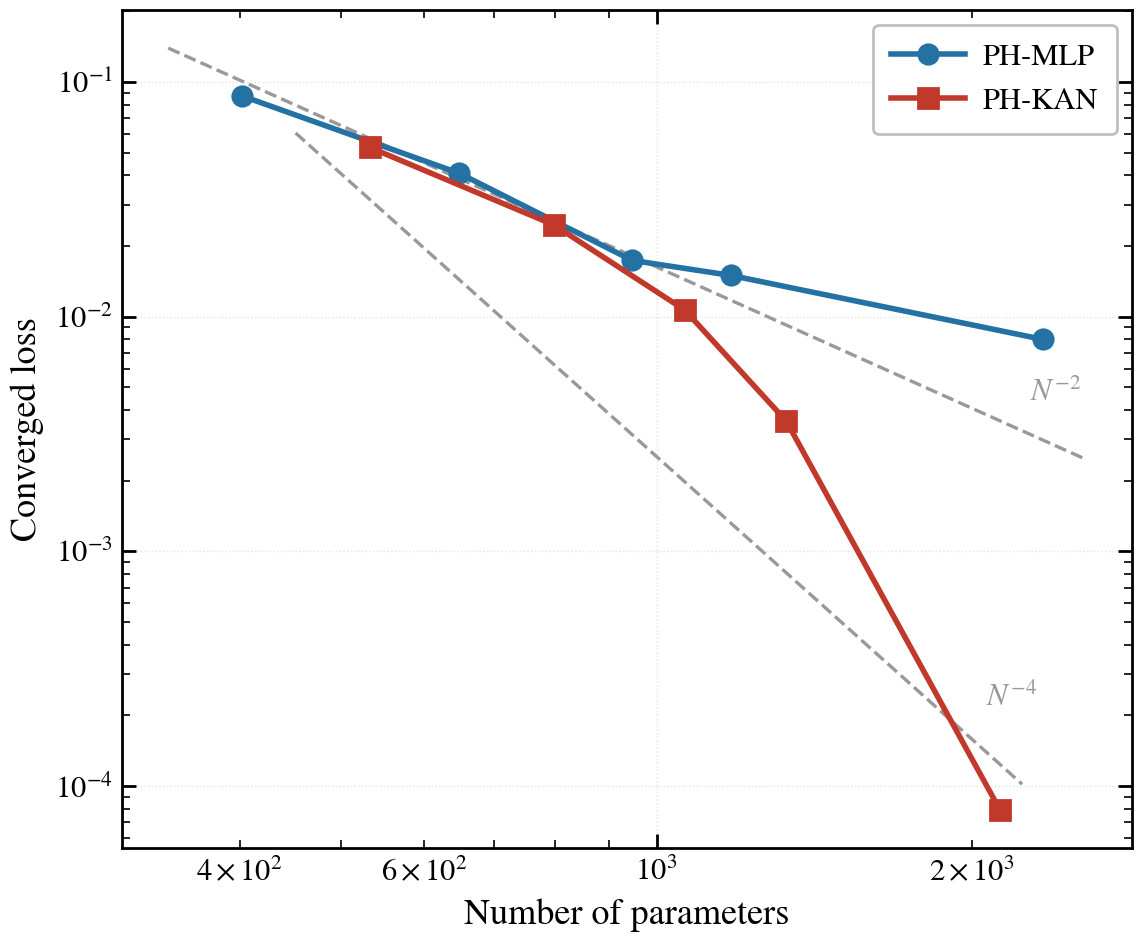}
    \caption{Test MSE as a function of the number of trainable parameters for PH-MLP and PH-KAN. The reference slopes $N^{-2}$ and $N^{-4}$ are shown for comparison.}
    \label{fig:scaling_comparison}
\end{figure}

As shown in Fig.~\ref{fig:scaling_comparison}, both models improve as the number of parameters increases, but the PH-KAN architecture exhibits a faster decay of the converged loss than the PH-MLP baseline. This is consistent with the favorable scaling behavior reported in the KAN literature, while showing that a similar trend remains visible in the present structure-preserving port-Hamiltonian identification setting.

\subsection{Interpretability Analysis}

One of the main motivations for using KANs instead of standard MLP parameterizations is that the learned functional structure can be directly inspected and simplified. Importantly, in the present paper we consider the strictly data-driven setting in which no prior structural information is assumed on the analytical forms of $J$, $R$, $H$, or $B$, beyond the port-Hamiltonian constraints enforced by construction. Even in this fully nonparametric setting, the learned KAN representations exhibit meaningful and inspectable internal structures after sparsification and pruning.

\begin{figure}[t]
    \centering
    \includegraphics[width=\linewidth]{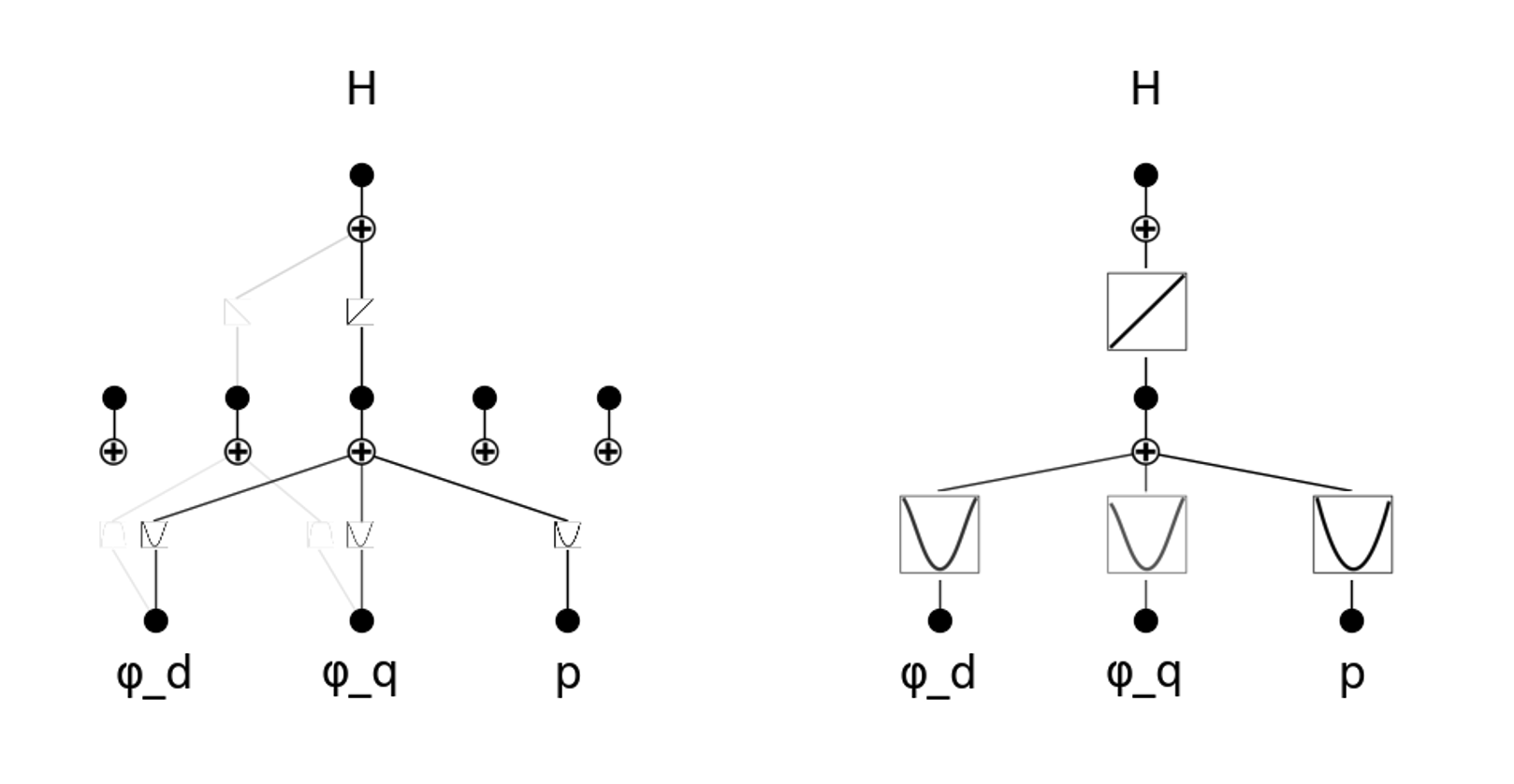}
    \caption{Learned KAN representation of the Hamiltonian. Left: model trained with sparsification regularization. Right: corresponding pruned representation.}
    \label{fig:hamiltonian_interpretability}
\end{figure}

\begin{figure}[t]
    \centering
    \includegraphics[width=\linewidth]{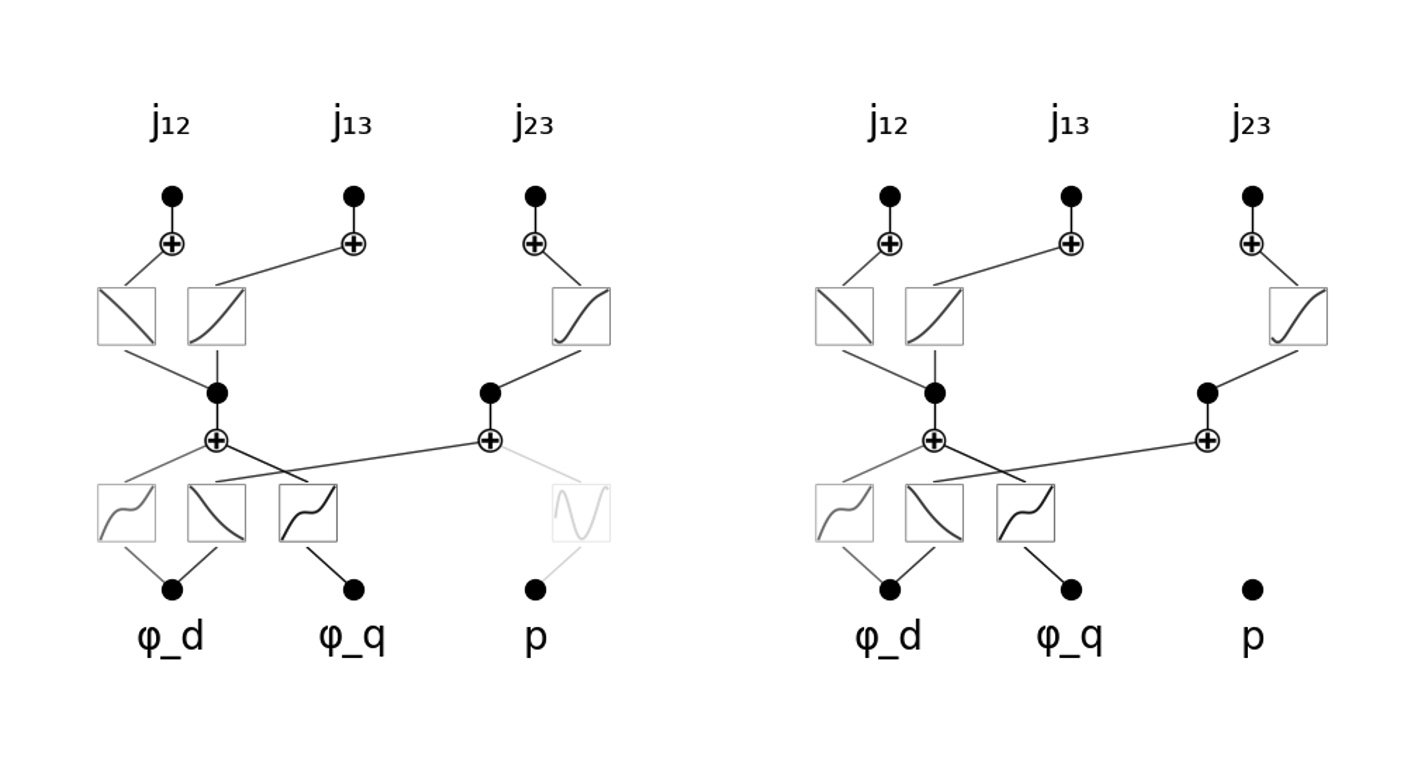}
    \caption{Learned KAN representation of the $J$ matrix. Left: model trained with sparsification regularization. Right: corresponding pruned representation.}
    \label{fig:j_interpretability}
\end{figure}

For the PMSM benchmark, the true Hamiltonian given in \eqref{eq:pmsm_matrices} is additively separable in the three state variables and quadratic in each of them. This structure is clearly reflected in the pruned KAN representation shown in Fig.~\ref{fig:hamiltonian_interpretability}: the dominant part of the learned map is organized through three main input-dependent branches associated with $\phi_d$, $\phi_q$, and $p$, followed by an aggregation and a final output transformation. After pruning, the resulting structure is consistent with a sum of one-dimensional contributions, each with a qualitatively quadratic behavior, which matches well the form of the true Hamiltonian.

A similar qualitative inspection can also be carried out for the interconnection matrix $J$. For the PMSM, the true matrix is state-dependent, but depends only on the electrical fluxes $\phi_d$ and $\phi_q$, and not on the momentum $p$; its nonzero entries are therefore determined by signed linear combinations of these electrical states. This key structural property is already reflected in the pruned representation of Fig.~\ref{fig:j_interpretability}, where the active graph is supported only by $\phi_d$ and $\phi_q$, while $p$ is completely inactive. Moreover, the two flux variables are first processed through similarly shaped, nearly linear edge functions, then combined into a shared intermediate representation, which is subsequently mapped to the outputs through edge functions with approximately opposite signs. Although this does not exactly recover the analytical expression of $J$ component by component, it remains qualitatively consistent with the signed linear state dependence of the true interconnection structure.

The fact that the exact analytical form of $J$ is not fully recovered is not surprising in the present setting, since the model is asked to learn $J$, $R$, $H$, and $B$ simultaneously without any additional prior knowledge on their individual forms. In such a fully flexible parameterization, part of the modeling burden can be redistributed across the different constitutive blocks through compensating nonlinear effects. A more faithful recovery of the true $J$ structure could likely be obtained if stronger prior information were incorporated, for instance by assuming a known dissipation or input matrix, or by constraining the admissible structure of $J$ itself and learning only its unknown coefficients.

A similar limitation applies to the dissipation matrix: in the present model, $R$ is learned as a state-dependent map, whereas the true PMSM dissipation matrix is constant. Its learned representation is therefore less directly interpretable from a physical viewpoint in this setting.

Overall, this illustrates an important advantage of KAN parameterizations over standard MLPs. In an MLP, the learned dependence is typically distributed across dense weight matrices and hidden activations, which makes such qualitative inspection difficult. By contrast, the KAN representation, especially after pruning, exposes the dominant functional organization of the learned constitutive maps and therefore provides a meaningful level of interpretability even in the absence of prior structural knowledge.

\section{Conclusion}
We presented a Port-Hamiltonian Kolmogorov--Arnold approach for structure-preserving identification of nonlinear port-Hamiltonian systems. The proposed PH-KAN model combines physical consistency with an interpretable learned representation. On the PMSM benchmark, it provides accurate trajectory prediction and yields constitutive structures whose dominant dependencies can be qualitatively inspected, with a Hamiltonian representation consistent with the true energy function. These results highlight the potential of KANs for interpretable nonlinear system identification. Future work will investigate how prior physical knowledge can be incorporated more explicitly into the PH-KAN architecture. In particular, when part of a constitutive relation is already known analytically, the corresponding branch of the network could be prescribed in advance, and only the remaining unknown coefficients or terms would be learned from data. A second research direction will consist of applying a similar approach using an irreversible port-Hamiltonian framework \cite{IPHS} to systems primarily governed by temperature and thermal effects.





\end{document}